\documentclass[11pt]{article}
\usepackage{amsmath}
\usepackage{amssymb}
\usepackage{graphicx, enumerate,verbatim}

\usepackage{theorem}
\numberwithin{equation}{section}

\newcommand{\RR} {\mathbb R}
\newcommand{\CC} {\mathbb C}
\newcommand{\ZZ} {\mathbb Z}
\newcommand{\NN} {\mathbb N}
\newcommand{\TT} {\mathbb T}
\newcommand{\DD} {\mathbb D}

\newcommand{\hC} {\widehat{\CC}}
\newcommand{\supp} {\operatorname{supp}}

\newcommand{\pa} {\partial}
\newcommand{\Cal} {\mathcal}

\newcommand{\beq} {\begin{equation}}
\newcommand{\eeq} {\end{equation}}
\newcommand{\demo} {\noindent {\it Proof. }}
\newcommand{\qed} {\hfill$\Box$ {\newline $\text{}$}}

\theoremstyle{plain}
\newtheorem{theorem}{Theorem}[section]
\newtheorem{proposition}[theorem]{Proposition}
\newtheorem{corollary}[theorem]{Corollary}
\newtheorem{lemma}[theorem]{Lemma}

\title{{The Maximal Entropy Measure of Fatou Boundaries}
\thanks{The second author was  partially supported by
NSF grant DMS-1500817
\newline {\it 2010 Math Subject Classification} 37F10, 30D05, 37A05
\newline
{\it Key words.} complex dynamics, Julia sets, Fatou sets, measure
 }}

\author{Jane Hawkins \\
Department of Mathematics, University of North Carolina \\
Chapel Hill, NC 27599, USA
\and Michael Taylor \\
Department of Mathematics, University of North Carolina \\
Chapel Hill, NC 27599, USA}

\date{}

\begin{document}

\maketitle

\begin{abstract}
We look at the maximal entropy (MME) measure of the boundaries of connected
components of the Fatou set of a rational map of degree $\ge 2$.  We show that
if there are infinitely many Fatou components, and if either the Julia set is
disconnected or the map is hyperbolic, then there can be at most one Fatou
component whose boundary has positive MME measure.
We also replace hyperbolicity by the more general hypothesis of
geometric finiteness.
\end{abstract}

\section{Introduction}\label{s1}

Let $R:\hC\rightarrow\hC$ be a rational map of degree $d \geq 2$,
defined on the Riemann sphere $\hC$, $\Cal{F}$ its
Fatou set, $\Cal{J}$ its Julia set, and $\lambda$ the unique maximal entropy
measure (MME) on $\Cal{J}$ (i.e., $h_\lambda(R)= \log d$).
Recall that $\Cal{F}$ and $\Cal{J}$ are invariant under $R$,
\beq
\text{$R$ preserves $\lambda$ and acts ergodically on }
(\Cal{J},\Cal{B},\lambda),
\label{1.1}
\eeq
where $\Cal{B}$ consists of the Borel subsets of $\Cal{J}$, and
\beq
\supp \lambda=\Cal{J}.
\label{1.2}
\eeq
See \cite{Lyu}, \cite{FLM}, or \S{5.4} of \cite{MNTU} for
\eqref{1.1}--\eqref{1.2}.
We want to study $\lambda(\pa\Cal{O})$, when $\Cal{O}$ is a connected component
of $\Cal{F}$, and in particular to understand when we can say $\lambda
(\pa\Cal{O})=0$.

The results in this paper hold for any Borel probability measure $\mu$
satisfying (\ref{1.1}) and (\ref{1.2})
It is well known that for many rational maps there are measures
in addition to the MME measure $\lambda$ that have these properties.
We focus on $\lambda$ to complement an earlier result by the authors
on the construction of the MME.  In the jointly written appendix of
\cite{Haw}, the authors showed
that the backward random iteration method works for drawing the Julia set
of any rational map. They proved the result by showing
the associated delta mass measures, averaged along almost
any randomly chosen backward path, converge to $\lambda$.
It is observed that certain details of Julia sets have a fuzzy appearance
for many maps when this algorithm is used, and the results of this paper
provide an explanation of this observation.

To get started, given a rational map $R$ and a component
$\Cal{O} \subset \Cal{F}$, Sullivan's Non-wandering Theorem
implies there exist $m\in\NN$ and a collection $\Cal{O}_1,\dots,\Cal{O}_k$
of mutually disjoint connected components of $\Cal{F}$ such that
\beq
R:\Cal{O}_j\longrightarrow \Cal{O}_{j+1},\quad j,j+1\in\ZZ/(k),
\label{1.3}
\eeq
and
\beq
R^m(\Cal{O})\subset \Cal{O}_1,\quad \text{hence }\ R^m(\pa\Cal{O})\subset
\pa\Cal{O}_1.
\label{1.4}
\eeq
Our first goal is to establish the following.

\begin{theorem} \label{t1.1}
Let $\Cal{O}_1,\dots,\Cal{O}_k$ be mutually disjoint
connected components of $\Cal{F}$ such that \eqref{1.3} holds.  Then either
\beq
\lambda(\pa\Cal{O}_j)=0,\quad \forall\, j\in\ZZ/(k),
\label{1.5}
\eeq
or
\beq
\pa\Cal{O}_1=\cdots=\pa\Cal{O}_k=\Cal{J}.
\label{1.6}
\eeq
If \eqref{1.5} holds, then $\lambda(\pa\Cal{O})=0$ whenever $\Cal{O}$ is a
component of $\Cal{F}$ such that $\Cal{O}\subset R^{-M}(\Cal{O}_j)$, $M \in \mathbb{N}$, and any $j$ .
\end{theorem}

The following general result will be useful in the proof of Theorem \ref{t1.1}.

\begin{lemma} \label{l1.2}
Let $(\Cal{J},\Cal{B},\lambda)$ be a probability space,
and assume $R:\Cal{J}\rightarrow \Cal{J}$ is $\Cal{B}$-measurable and
$\lambda$ is ergodic and invariant under $R$. Then, for $E\in\Cal{B}$,
\beq
R^{-1}(E)\supset E\Longrightarrow \lambda(E)=0\ \text{ or }\ 1.
\label{1.7}
\eeq
\end{lemma}

\demo
Let $F_j=R^{-j}(E)$.  Then $E=F_0\subset F_j\nearrow F\in\Cal{B}$.
The hypotheses also imply $\lambda(F)=\lambda(E)$ and $R^{-1}(F)=F$.
Ergodicity implies $\lambda(F)=0$ or $1$, so we have \eqref{1.7}.
\qed

\noindent {\it Proof of Theorem \ref{t1.1}.\ }
Lemma \ref{l1.2} applies to $E=\pa\Cal{O}_1\cup\cdots\cup\pa\Cal{O}_k$, and we
have \eqref{1.7}.  Also $R^{-1}(\pa\Cal{O}_{j+1})\supset\pa\Cal{O}_j$, so
\beq
\lambda(\pa\Cal{O}_{j+1})\ge \lambda(\pa\Cal{O}_j),\quad \forall\, j\in\ZZ/(k).
\label{1.8}
\eeq
This implies
\beq
\lambda(\pa\Cal{O}_1)=\cdots=\lambda(\pa\Cal{O}_k).
\label{1.9}
\eeq
Applying \eqref{1.7} leads to two cases:
\beq
\lambda(\pa\Cal{O}_1\cup\cdots\cup\pa\Cal{O}_k)=0,\ \text{ or }\
\lambda(\pa\Cal{O}_1\cup\cdots\cup\pa\Cal{O}_k)=1.
\label{1.10}
\eeq
In the first case, \eqref{1.5} holds.  Furthermore,
$\lambda(\pa\Cal{O})=0$ for each $\Cal{O}\subset R^{-M}
(\Cal{O}_j)$, $M\in\NN$, since
\beq
\aligned
R^M(\Cal{O})=\Cal{O}_j&\Rightarrow
R^M(\overline{\Cal{O}})=\overline{\Cal{O}}_j\ \text{and}\
R^M(\pa\Cal{O})\subset \pa\Cal{O}_j \\
&\Rightarrow \pa\Cal{O}\subset R^{-M}(\pa\Cal{O}_j).
\endaligned
\label{1.11}
\eeq
In the second case, since $\pa\Cal{O}_1\cup\cdots\cup
\pa\Cal{O}_k$ is compact and \eqref{1.2} holds, we have
\beq
\pa\Cal{O}_1\cup\cdots\cup\pa\Cal{O}_k=\Cal{J}.
\label{1.12}
\eeq

Going further, let us note that $R$ and $R^k$ have the same Julia set and the
same maximal entropy measure.  Given (1.3), we have $R^k:\Cal{O}_j\rightarrow
\Cal{O}_j$ for each $j\in\ZZ/(k)$, hence $R^k:\pa\Cal{O}_j\rightarrow
\pa\Cal{O}_j$.  Hence, by Lemma 1.2 (applied to $R^k$), for {\it each} such $j$,
\beq
\lambda(\pa\Cal{O}_j)=0\ \text{or}\ 1.
\label{1.13}
\eeq
Then
\beq
\lambda(\pa\Cal{O}_j)=1\Longrightarrow \pa\Cal{O}_j=\Cal{J}.
\label{1.14}
\eeq
We still have \eqref{1.9}, so
\beq
\lambda(\pa\Cal{O}_j)=1\Longrightarrow \pa\Cal{O}_1=\cdots =\pa\Cal{O}_k
=\Cal{J}.
\label{1.15}
\eeq
This proves Theorem \ref{t1.1}.

\qed

The following result of \cite{CMT} is a significant consequence of Theorem
\ref{t1.1}.

\begin{theorem} \label{c1.3}
Either
\beq
\lambda(\pa\Cal{O})=0\ \text{for each connected component $\Cal{O}$ of }
\Cal{F},
\label{1.16}
\eeq
or there is a connected component $\Cal{O}_1$ of $\Cal{F}$ such that
\beq
\pa\Cal{O}_1=\Cal{J}.
\label{1.17}
\eeq
More precisely, $\Cal{O}_1$ can be taken as in \eqref{1.3}--\eqref{1.4}, and
then \eqref{1.6} holds.
\end{theorem}

We can rephrase Theorem \ref{c1.3} using the notion of {\it residual Julia set},
defined by
\beq
\Cal{J}_0=\Cal{J}\setminus \bigcup\limits_j \pa\Cal{O}_j,
\label{1.18}
\eeq
where $\{\Cal{O}_j\}$ consists of all the connected components of $\Cal{F}$.
Note that
\beq
\eqref{1.16}\Rightarrow \Cal{J}_0\neq\emptyset,\quad \text{and }\
\eqref{1.17}\Rightarrow \Cal{J}_0=\emptyset.
\label{1.19}
\eeq
As noted in \cite{CMT}, it follows from Theorem \ref{c1.3} that
\beq
\Cal{J}_0\neq\emptyset\Longrightarrow \lambda(\Cal{J}_0)=1,
\label{1.20}
\eeq
and we have the following basic result of \cite{Mor}.

\begin{corollary} \label{c1.4}
Either $\Cal{J}_0\neq\emptyset$
or there is a Fatou component $\Cal{O}_1$ such that $\pa
\Cal{O}_1=\Cal{J}$.  If $\Cal{J}_0\neq\emptyset$, then $\Cal{J}_0$ is a
dense, $\Cal{G}_\delta$ subset of $\Cal{J}$.
\end{corollary}

To establish the results in the last sentence of Corollary \ref{c1.4}, we note that
denseness of $\Cal{J}_0$ in $\Cal{J}$ follows
from \eqref{1.20} and \eqref{1.2}.  The fact
that $\Cal{J}_0$ is obtained from $\Cal{J}$ by successively removing
$\pa\Cal{O}_j$, $j\in\NN$, implies $\Cal{J}_0$ is a $\Cal{G}_\delta$ subset
of $\Cal{J}$.

$\text{}$

Our goal is the study of $\lambda(\pa\Cal{O})$ for various Fatou
components in cases where this measure is not identically zero on these
boundaries.  An important class of rational maps with empty residual
Julia set is the class of polynomials of degree $d\ge 2$, which we tackle in
\S{\ref{s2}}.  In this case, the Fatou component $\Cal{O}^\infty$ containing $\infty$
satisfies
\beq
\pa\Cal{O}^\infty=\Cal{J}.
\label{1.21}
\eeq
We give conditions under which we can show that,
if $\Cal{O}$ is another component of
$\Cal{F}$, i.e., a bounded component of $\Cal{F}$, then
\beq
\lambda(\pa\Cal{O})=0.
\label{1.22}
\eeq
We demonstrate this for a class of polynomials whose Fatou sets have an infinite
number of components.  One important property of $\Cal{O}^\infty$ used in the
analysis is its complete invariance: $R(\Cal{O}^\infty)\subset\Cal{O}^\infty$
and $R^{-1}(\Cal{O}^\infty)\subset\Cal{O}^\infty$.

In \S{\ref{s3}} we extend the scope of this work to include other rational maps $R$
for which there is a completely invariant Fatou component $\Cal{O}^\#$, so
\beq
R(\Cal{O}^\#)\subset\Cal{O}^\#\ \text{ and }\ R^{-1}(\Cal{O}^\#)\subset
\Cal{O}^\#.
\label{1.23}
\eeq
We extend the basic results of \S{\ref{s2}} to cover this more general
situation.  In \S{\ref{s4}} we consider hyperbolic maps and
geometrically finite maps.
Results of \S\S{\ref{s3}--\ref{s4}} together yield the following.

\begin{theorem} \label{t1.5}
Let $R$ be a rational map of degree $\ge 2$, and
assume $\Cal{F}$ has infinitely many connected components.  Assume
there is a Fatou component $\Cal{O}_1$ such that $\lambda(\pa\Cal{O}_1)\neq
0$.  Then $\lambda(\pa\Cal{O})=0$ for each Fatou component $\Cal{O}\neq
\Cal{O}_1$, under either of the following conditions:
\beq
\text{$\Cal{J}$ is disconnected,}
\tag{A}
\eeq
or
\beq
\text{$\Cal{J}$ is connected and $R$ is hyperbolic.}
\tag{B}
\eeq
Moreover, if (A) or (B) hold, then $\pa\Cal{O}_1=\Cal{J}$, and $\Cal{O}_1$
is completely invariant.   
\end{theorem}

Theorem \ref{t1.5} follows from Propositions \ref{p3.3} and \ref{p4.1}.
We mention that Proposition \ref{p3.2} has the same conclusion, with
hypothesis (B) replaced by
\beq
\aligned
&\text{$\Cal{J}$ is connected and locally connected, and} \\
&\text{$\Cal{F}$ has a completely invariant component.}
\endaligned
\tag{B1}
\eeq
The special case when $R$ is a polynomial and $\Cal{J}$ is connected and
locally connected is Proposition \ref{p2.3}.
As we show in \S{\ref{s4}}, using (B1), we can also replace (B) by
\beq
\text{$\Cal{J}$ is connected and $R$ is geometrically finite.}
\tag{B2}
\eeq

We end with an appendix, giving examples to illustrate our results.

\section{Polynomials} \label{s2}

Here we explore consequences of Theorem \ref{t1.1} for the class
of polynomials of degree $d\ge 2$, e.g.,
\beq
R(z)=z^d+a_{d-1}z^{d-1}+\cdots+a_0,\quad d\ge 2.
\label{2.1}
\eeq
(If the leading term were $a_dz^d$, we could conjugate by $z\mapsto cz$ with
$c^{d-1}=a_d$, to obtain the form \eqref{2.1}.)  In such a case, $\Cal{F}$
has a connected component $\Cal{O}^\infty$, containing $\infty$, and
(cf.~\cite{B}, Theorem 5.2.1)
\beq
\Cal{J}=\pa\Cal{O}^\infty.
\label{2.2}
\eeq
It clearly follows from (\ref{2.2}) that for polynomials,
$$\Cal{J}_0 = \emptyset \quad
\rm{and} \quad \lambda(\pa\Cal{O}^\infty) = 1.
$$
The set $\Cal{K}=\hC\setminus\Cal{O}^\infty$ is called the
{\it filled Julia set} of $R$, and one has $\Cal{J}=\pa\Cal{K}$.
In some cases, $\Cal{F}=\Cal{O}^\infty$.
If $\Cal{F}$ has exactly 2 connected components, say $\Cal{O}^\infty$ and
$\Cal{O}_1$, then it is also the case that $\pa\Cal{O}_1=\Cal{J}$.  Otherwise,
$\Cal{F}$ must have an infinite number of connected components
(cf.~\cite{B}, Theorem 5.6.2).  We call the
connected components of $\Cal{F}$ other than $\Cal{O}^\infty$, i.e., those
contained in $\Cal{K}$, {\it bounded Fatou components}.  The component
$\Cal{O}^\infty$ is completely invariant.  Hence, if $\Cal{O}$
is a bounded Fatou component, so is $R(\Cal{O})$ and so is each connected
component of $R^{-1}(\Cal{O})$.  We now look into when a bounded Fatou
component $\Cal{O}$ can be shown to satisfy $\lambda(\pa\Cal{O})=0$.

\begin{proposition} \label{p2.1}
Let $R(z)$ have the form \eqref{2.1}, and filled Julia
set $\Cal{K}$.  Assume there is a point $z_0\in\Cal{K}$ such that $\Cal{K}
\setminus\{z_0\}$ is not connected.
Then $\lambda(\pa\Cal{O})=0$ for each
bounded Fatou component $\Cal{O}$ of $R$.
\end{proposition}

\demo
Take a bounded Fatou component $\Cal{O}\subset\Cal{K}$.
The non-wandering theorem gives $m\in\NN$ and Fatou components $\Cal{O}_1,
\dots,\Cal{O}_k$ (necessarily bounded Fatou components)
such that \eqref{1.3}--\eqref{1.4} hold.
Say $\Cal{K}_1$ and $\Cal{K}_2$ are disjoint connected components of
$\Cal{K}\setminus\{z_0\}$.
We can assume $\Cal{O}_1\subset\Cal{K}_1$.  Then $\overline{\Cal{O}}_1$
does not contain $\pa\Cal{K}_2$, so \eqref{1.6} fails.
By Theorem \ref{t1.1}, this forces
$\lambda(\pa\Cal{O}_1)=0$, and hence $\lambda(\pa\Cal{O})=0$.
\qed

Of course, Proposition 2.1 applies if $\Cal{K}$ is not connected.  We will now
assume $\Cal{K}$ is connected.  Hence $\hC\setminus\Cal{K}$ is simply connected.
By the Riemann mapping theorem, there is a unique biholomorphic map
\beq
\varphi:\hC\setminus \overline{\DD}\longrightarrow \hC\setminus\Cal{K},
\label{2.3}
\eeq
such that
\beq
\varphi(\infty)=\infty,\quad D\varphi(\infty)=\alpha>0.
\label{2.4}
\eeq
Here $\overline{\DD}$ is the disk $\{z\in\CC:|z|\le 1\}$.
We use this to recall
from Chapter 6 of \cite{McM} conditions under which
the hypothesis of Proposition \ref{p2.1} holds.
It involves the notion of external rays,
\beq
\psi_\theta(r)=\varphi(r e^{i\theta}),\quad r\in (1,\infty),\
\theta\in\TT^1=\RR/(2\pi \ZZ).
\label{2.5}
\eeq
Given $\theta\in\TT^1$, we say $\psi_\theta$ lands at a point $z\in\pa\Cal{K}$
provided
\beq
\lim\limits_{r\searrow 1}\, \psi_\theta(r)=z.
\label{2.6}
\eeq
We mention two classical results, given as Theorems 6.1 and 6.2 of \cite{McM}.
The first is that the limit in \eqref{2.6}
exists for almost all $\theta\in\TT^1$.
The second is the following:
\beq
\aligned
&\text{If $E\subset \TT^1$ has positive measure, then there exist
$\theta_0,\theta_1\in E$ such that} \\
&\lim\limits_{r\searrow 1}\, \psi_{\theta_0}(r)=z_0,\ \
\lim\limits_{r\searrow 1}\, \psi_{\theta_1}(r)=z_1,\ \ \text{and }\
z_0\neq z_1.
\endaligned
\label{2.7}
\eeq
These results lead to the following, which is part of Corollary 6.7 of
\cite{McM}.

\begin{proposition} \label{p2.2}
Let $R$ have the form \eqref{2.1} and assume $\Cal{K}$ is
connected.  Take $\varphi$ and $\psi_\theta$ as in \eqref{2.3}--\eqref{2.5}.
Assume there exist $\xi_0\neq \xi_1\in\TT^1$ such that
\beq
\lim\limits_{r\searrow 1}\, \psi_{\xi_0}(r)=
\lim\limits_{r\searrow 1}\, \psi_{\xi_1}(r)=z_0.
\label{2.8}
\eeq
Then $\Cal{K}\setminus \{z_0\}$ is not connected.
\end{proposition}

\demo
The union $\gamma$ of the two rays
$\{\psi_{\xi_0}(r):1<r<\infty\}$ and $\{\psi_{\xi_1}(r):1<r<\infty\}$, together
with their endpoints $z_0$ and $\infty$, forms a simple closed curve in $\hC$.
The Jordan curve theorem implies that $\hC\setminus \gamma$ has exactly two
connected components.  Our hypotheses imply
\beq
\Cal{K}\cap \gamma=\{z_0\}.
\label{2.9}
\eeq
It remains to note that each connected component of $\hC\setminus\gamma$
contains a point of $\Cal{K}$, and this follows from \eqref{2.7}, first taking $E$
to be the open arc from $\xi_0$ to $\xi_1$ in $\TT^1$, then taking $E$ to be
the complementary open arc in $\TT^1$.
\qed

Regarding the question of when Proposition \ref{p2.2} applies, we note that
a definite answer can be given under the additional hypothesis that
\beq
\text{$\Cal{J}$ is connected and locally connected,}
\label{2.10}
\eeq
or equivalently, that $\Cal{K}$ is connected and locally connected.  In that
case, a classical result of Caratheodory (cf.~\cite{Mil}, Theorem 17.14) implies
that $\varphi$ in \eqref{2.3}, mapping
$\DD^\infty=\hC\setminus{\overline{\DD}}$ to
$\Cal{O}^\infty=\hC\setminus \Cal{K}$, has a unique continuous extension to
\beq
\varphi:\overline{\DD}^\infty\longrightarrow \overline{\Cal{O}}^\infty.
\label{2.11}
\eeq
The image of $\varphi$ in \eqref{2.11} is both compact and dense in
$\overline{\Cal{O}}^\infty$, so $\varphi$ in \eqref{2.11} is surjective.  It
restricts to a continuous map
\beq
\varphi:\pa \DD\longrightarrow \Cal{J},
\label{2.12}
\eeq
also surjective.  In view of \eqref{2.5}, we have
\beq
\varphi(e^{i\theta})=\lim\limits_{r\searrow 1}\, \psi_\theta(r).
\label{2.13}
\eeq
Now if $\varphi$ in \eqref{2.12} is also one-to-one, this makes $\Cal{J}\subset
\hC$ a simple closed curve, so $\Cal{F}=\hC\setminus\Cal{J}$ would have
just two connected components.  We have the following conclusion.

\begin{proposition} \label{p2.3}
Let $R$ have the form \eqref{2.1} and assume $\Cal{J}$
is connected and locally connected.  If $\Cal{F}$ has infinitely many
connected components, then
\beq
\lambda(\pa\Cal{O})=0\ \text{ for each bounded component $\Cal{O}$
of $\Cal{F}$}.
\label{2.14}
\eeq
\end{proposition}

\demo
As we have just seen, the hypotheses imply that $\varphi$ in
\eqref{2.12} is not one-to-one.  Hence there exist $\xi_0\neq\xi_1\in\TT^1$
such that $\varphi(e^{i\xi_0})=\varphi(e^{i\xi_1})$ in \eqref{2.13}.
Thus, with $z_0=\varphi(e^{i\xi_0})=\varphi(e^{i\xi_1})$,
Proposition \ref{p3.2} implies that $\Cal{K}\setminus\{z_0\}$
is not connected, so \eqref{2.14} follows from Proposition \ref{p2.1}.
\qed

\section{Other maps with a completely invariant Fatou component}\label{s3}

Extending our scope a bit, let us now assume that $R$ (of degree $d\ge 2$)
has the property that $\Cal{F}$ has a completely invariant connected
component $\Cal{O}^\#$, i.e.,
\beq
R(\Cal{O}^\#)\subset\Cal{O}^\#\ \text{ and }\ R^{-1}(\Cal{O}^\#)\subset
\Cal{O}^\#.
\label{3.1}
\eeq
If there is a fixed point $p\in\Cal{O}^\#$ of $R$, then, conjugating by a
linear fractional transformation, we can take $p=\infty$.
If  $R^{-1}(p)=p$ as well, then $R$ is a polynomial
(which guarantees that $p$ must be a superattracting fixed point).
There are many examples of rational maps with completely invariant Fatou
domains and attracting fixed points that are not polynomials.
We discuss some in Appendix \ref{sa}.

Generally when \eqref{3.1} holds, we have
\beq
\pa\Cal{O}^\#=\Cal{J}.
\label{3.2}
\eeq
This follows from Theorem 5.2.1 of \cite{B},
which also contains the results that
$\Cal{O}^\#$ is either simply connected or infinitely connected;
all the other components of $\Cal{F}$ are simply connected;
and $\Cal{O}^\#$ is simply
connected if and only if $\Cal{J}$ is connected.
Using these facts, we can set things up to obtain results parallel to
Propositions \ref{p2.1}--\ref{p2.3}.
To formulate these results, let us arrange that
$\infty\in\Cal{O}^\#$.  Then $\Cal{O}^\#$ is the unbounded component of
$\Cal{F}$, and we call the other components of $\Cal{F}$ bounded Fatou
components.  As in the case of polynomials, we see that if $\Cal{O}$ is a
bounded Fatou component, so is $R(\Cal{O})$ and so is each connected component
of $R^{-1}(\Cal{O})$.  We set $\Cal{K}=\hC\setminus\Cal{O}^\#$, and call this
the filled Julia set of $R$.  The following result has the same proof as
Proposition \ref{p2.1}.

\begin{proposition} \label{p3.1}
Let $R$ be a rational map of degree $\ge 2$ for which
a Fatou component $\Cal{O}^\#$ satisfies \eqref{3.1},
and let $\Cal{K}=\hC\setminus\Cal{O}^\#$ be its filled Julia set.
Assume there exists $z_0\in\Cal{K}$
such that $\Cal{K}\setminus\{z_0\}$ is not connected.  Then $\lambda(\pa\Cal{O})
=0$ for each Fatou component $\Cal{O}\neq\Cal{O}^\#$.
\end{proposition}

As before, Proposition \ref{p3.1} certainly applies if $\Cal{K}$ is not connected.
If $\Cal{K}$ is connected, then $\Cal{O}^\#$ is simply connected, and we
again have the set-up \eqref{2.3}--\eqref{2.7},
and Proposition \ref{p2.2} extends to this setting, as does the analysis
leading to the following extension of Proposition \ref{p2.3}.

\begin{proposition} \label{p3.2}
Let $R$ be a rational map of degree $\ge 2$ for which
there is a Fatou component $\Cal{O}^\#$ satisfying \eqref{3.1},
and assume $\Cal{J}$ is connected and locally connected.
If $\Cal{F}$ has infinitely many
connected components, then $\lambda(\pa\Cal{O})=0$ for each Fatou component
$\Cal{O}\neq \Cal{O}^\#$.
\end{proposition}

In counterpoint, we have the following result when $\Cal{J}$ is not connected.

\begin{proposition} \label{p3.3}
Let $R$ be a rational map of degree $\ge 2$, and
assume $\Cal{J}$ is disconnected.  If there is a Fatou component $\Cal{O}_1$
such that $\lambda(\pa\Cal{O}_1)\neq 0$,
then $\lambda(\pa\Cal{O})=0$ for each
Fatou component $\Cal{O}\neq\Cal{O}_1$.
\end{proposition}

\demo
The hypotheses imply that the residual set $\Cal{J}_0$ is empty.
Hence, given given $\Cal{J}$ disconnected,
by Theorem 4.4.19 of \cite{MNTU},
$\Cal{F}$ has a completely invariant component $\Cal{O}^\#$.
By Proposition \ref{p3.1}, $\lambda(\pa\Cal{O})=0$ for each Fatou
component $\Cal{O}\neq\Cal{O}^\#$.  This proves the proposition (and forces
$\Cal{O}_1=\Cal{O}^\#$).
\qed

$\text{}$ \newline
{\sc Remark.} As seen in \S{\ref{s1}}, the hypothesis that
$\lambda(\pa\Cal{O}_1)\neq 0$ for some Fatou component $\Cal{O}_1$ is
equivalent to the hypothesis that $\Cal{J}_0=\emptyset$.  The {\it Makienko
conjecture} can be stated as saying that, if $\Cal{J}_0=\emptyset$,
then there is a Fatou component that is completely invariant under $R^2$.
A discussion of conditions under which Makienko's conjecture has been
proved can be found in \cite{CMMR}.  Of particular use here is the following
result of \cite{Q}:
\beq
\text{The Makienko conjecture holds provided $\Cal{J}$ is locally connected.}
\label{3.3}
\eeq
We make use of this in the following section.

\section{Hyperbolic maps and geometrically finite maps} \label{s4}

A rational map $R:\hC\rightarrow \hC$, of degree $d\ge 2$, with Julia set
$\Cal{J}$, is said to be {\it hyperbolic} provided
\beq
\Cal{J}\cap \overline{C^+(R)}=\emptyset,\quad \text{where }\
C^+(R)=\bigcup\limits_{k\ge 1} R^k(C_R),
\label{4.1} \eeq
where $C_R=\{z\in\hC:DR(z)=0\}$ is the set of critical points of $R$.  An
equivalent condition for hyperbolicity (cf.~\cite{MNTU}, Theorem 4.4.2)
is that each critical point of $R$ is in $\Cal{F}$
and each forward orbit of a critical point converges
to an attracting cycle.  Hyperbolic rational maps are relatively ``tame,''
from a topological point of view.
For example (cf.~\cite{MNTU}, Theorem 4.4.5), if $R$ is hyperbolic,
\beq
\text{$\Cal{J}$ connected $\Longrightarrow \Cal{J}$ locally connected.}
\label{4.2}
\eeq
There is also a partial converse to the implication
\eqref{3.1} $\Rightarrow$ \eqref{3.2}
for hyperbolic maps (cf.~\cite{MNTU}, Theorem 4.4.16):
\beq
\gathered
\text{If $R$ is hyperbolic and $\Cal{J}$ is connected, and $\Cal{O}$
is a component of $\Cal{F}$, then} \\
R(\Cal{O})\subset\Cal{O},\ \pa\Cal{O}=\Cal{J}\Longrightarrow \Cal{O}\
\text{is completely invariant.}
\endgathered
\label{4.3}
\eeq
We therefore have the following.

\begin{proposition} \label{p4.1}
Assume that $R$ is hyperbolic and $\Cal{J}$ is
connected, that $\Cal{F}$ has infinitely many connected components,
and that there is a Fatou component $\Cal{O}_1$
such that $\lambda(\pa\Cal{O}_1) \neq 0$.
Then $\lambda(\pa\Cal{O})=0$ for each Fatou component $\Cal{O}\neq
\Cal{O}_1$.
\end{proposition}

\demo
By Theorem \ref{c1.3}, there must be a component $\Cal{O}^\#$ of $\Cal{F}$
such that $\pa\Cal{O}^\#=\Cal{J}$.  Furthermore, $\Cal{O}^\#$ is invariant
under $R^k$, which is also hyperbolic.
By \eqref{4.3}, such $\Cal{O}^\#$ must be
completely invariant (under $R^k$).
Thanks to \eqref{4.2}, Proposition \ref{p3.2}
applies (to $R^k$), and it implies
that $\lambda(\pa\Cal{O})=0$ for each component $\Cal{O}\neq\Cal{O}^\#$.

Incidentally, this forces $\Cal{O}_1=\Cal{O}^\#$.  Furthermore, taking into
account \eqref{1.6}, we see that $k$ must be $1$.
Hence $\Cal{O}_1$ must be completely invariant under $R$.
\qed

For further results, we bring in the following two generalizations of
hyperbolicity, taking definitions from \cite{MNTU}, p.~153.

$\text{}$ \newline
{\bf Definition.} A rational map $R:\hC\rightarrow\hC$ of degree $\ge 2$ is
said to be {\it subhyperbolic} provided the following two conditions hold:
\beq
\aligned
&\text{the forward orbit of each critical point in $\Cal{F}$} \\
&\text{converges to an attracting cycle},
\endaligned
\label{4.4}
\eeq
and
\beq
\aligned
&\text{the forward orbit of each critical point in $\Cal{J}$} \\
&\text{is eventually periodic.}
\endaligned
\label{4.5}                                                             
\eeq
If we merely assume \eqref{4.5} holds, we say $R$ is {\it geometrically finite}.

$\text{}$

Clearly $R$ hyperbolic $\Rightarrow R$ subhyperbolic $\Rightarrow R$
geometrically finite.  The usefulness of geometric finiteness for the work
here stems from the following generalization of \eqref{4.2}, established
in \cite{Tan}, Theorem A:
\beq
\aligned
&\text{If $R$ is geometrically finite and $\Cal{J}$ is connected,} \\
&\text{then $\Cal{J}$ is locally connected.}
\endaligned
\label{4.6}
\eeq
In concert with \eqref{3.3}, this yields the following:

\begin{lemma} \label{l4.2}
Assume that $R$ is geometrically finite and $\Cal{J}$ is connected.  Then
either $\Cal{J}_0\neq \emptyset$ or $\Cal{F}$ has a completely invariant
component for $R^2$.
\end{lemma}

We now establish the following extension of Proposition \ref{p4.1}.

\begin{proposition} \label{p4.3}
Assume $R$ is geometrically finite, $\Cal{J}$ is connected, and $\Cal{F}$
has infinitely many connected components.
Assume further that there is a Fatou component
$\Cal{O}_1$ such that $\lambda(\pa\Cal{O}_1)\neq 0$.  Then $\lambda(\pa\Cal{O})
=0$ for each Fatou component $\Cal{O}\neq\Cal{O}_1$.
\end{proposition}

\demo
The hypotheses imply $\Cal{J}_0=\emptyset$, by \eqref{1.20}.  Hence, by
Lemma \ref{l4.2}, there is a component $\Cal{O}^\#$ of $\Cal{F}$ that is
completely invariant under $R^2$.  As in \eqref{3.2}, this implies $\pa
\Cal{O}^\#=\Cal{J}$.  We are assuming $\Cal{J}$ is connected, and by
\eqref{4.6} $\Cal{J}$ is also locally connected.  Thus Proposition \ref{p3.2}
applies (to $R^2$), giving $\lambda(\pa\Cal{O})=0$ for each Fatou component
$\Cal{O}\neq\Cal{O}^\#$.

As before, this forces $\Cal{O}_1=\Cal{O}^\#$.  It also implies that
$\Cal{O}^\#$ is completely invariant under $R$.
\qed

$\text{}$

We turn to a proof of Theorem 1.5.

\noindent {\it Proof of Theorem 1.5}.     
If (A) holds, the conclusion that $\lambda(\pa\Cal{O})=0$ for each Fatou
component $\Cal{O}\neq\Cal{O}_1$ is a direct consequence of Proposition
\ref{p3.3}.  The proof of Proposition \ref{p3.3} also ends with the comment
that $\Cal{O}_1$ must coincide with a completely invariant component
$\Cal{O}^\#$.  Similarly, if (B) holds, the conclusion of the theorem follows
from Proposition \ref{p4.1}, the complete invariance of $\Cal{O}_1$
again established in the course of the proof of Proposition \ref{p4.1}.
Validity of Theorem \ref{t1.5} with (B) replaced by (B1) follows from
Proposition \ref{p3.2}, and its validity with (B) replaced by (B2) follows from
Proposition \ref{p4.3}.
\qed

\appendix
\section{Examples} \label{sa}

We include a few basic examples illustrating the results of
\S\S{\ref{s1}--\ref{s4}}.  In all these examples, $\Cal{J}_0=\emptyset$.
Before getting to them, we make a few useful general comments whose proofs are
found in the literature in the bibliography.  If $R$ is a rational map of degree $d \geq 2$ with 
a $k$-periodic point $z_0$ whose multiplier satisfies:  $DR^k(z_0)=\exp(2 \pi i / m)$,
$m \in \mathbb{N}$,  then the $k$-cycle is called {\em parabolic}; $k$ always refers to the minimum
period.

$\text{}$ \newline
{\bf Remarks A.} 

\begin{enumerate} 

\item  If a rational map $R$  has an attracting $k$-periodic point in $\Cal{F}$, then there are at least
$k$ Fatou components in $\Cal{F}$.

\item If $R$ has a parabolic $k$-periodic point  $z_0$ with multiplier  $DR^k(z_0)=\exp(2 \pi i / m)$,
then there are at least $mk$ Fatou components in $\Cal{F}$. 

 \item  These lead to the following lemma.
\begin{lemma} \label{la.1}
Let $R:\hC\rightarrow\hC$ be a rational map of degree $\ge 2$.  Assume
its Fatou set $\Cal{F}$ has a completely invariant component and that
$R$  has an attracting or parabolic  $k$-periodic point  with $k\ge 2$.  Then
$\Cal{F}$ has infinitely many components.
\end{lemma}

\demo
The completely invariant component cannot be one of the components in the $n$-cycle,
$n \geq 2$, of components in $\Cal{F}$ induced by the non-repelling $k$-cycle for $R$
so $\Cal{F}$ has more than two components.
\qed

\item Whenever $R:\hC\rightarrow\hC$ is a rational map of degree $2$,
the Julia set $\Cal{J}$ is either connected or totally disconnected.

\item If every critical point in $\Cal{J}$ has a finite forward orbit (i.e., $R$ is geometrically finite),  then all non-repelling periodic points of $R$ are attracting or parabolic. 
If $R$ has degree $d \geq 2$ and there are $2d-2$ critical points in $\Cal{F}$, then $R$ is geometrically finite. 

\end{enumerate}

$\text{}$ \newline
{\bf Example 1.}   As an example of a non-polynomial map with a completely invariant Fatou component, neither hyperbolic nor subhyperbolic, 
but geometrically finite, we set:
\beq \label{eqn:quad2}
R(z) =  z+\frac{1}{z} + \frac 3 2. \quad
\eeq
$R$ has critical points $c_1=1$ and $c_2=-1$. 

For these maps $\infty$ is a fixed point with multiplier $1$ so the map is parabolic; this implies that $\infty, 0 \in \Cal J$ and there is an attracting petal in $\Cal{F}$
which sits in a component $\Cal{O}^\# \subset \Cal{F}$ containing a critical point. We have:
$$-1 \mapsto -1/2 \mapsto -1$$  
so there is a superattracting 2-cycle of components in $\Cal{F}$, $\Cal{O}_1$ and $\Cal{O}_2= R(\Cal{O}_1)$.  Therefore $c_1 \in  \Cal{O}^\# $. 
Since there are no other critical points, $R$ is geometrically finite by Remark A.5.

There are infinitely many
components $\Cal{O} \subset \Cal{F}$; in fact all $\Cal{O} \neq  \Cal{O}^\#$  map to $\Cal{O}_1 \cup \Cal{O}_2$; 
the first statement follows from Lemma \ref{la.1}, and the second
since we have exhausted the critical points so no other Fatou components are possible.  Every component of $\Cal{F}$ is simply connected by
Remark A.4.

Since  $c_1 \in \Cal{O}^\#$ and $R(\Cal{O}^\#)=\Cal{O}^\#$,      $\Cal{O}^\#$ is a degree 2 branched cover of itself, so $\Cal{O}^\#$ is completely invariant.  
The map $R$ is geometrically finite, and $\Cal{J}$ is connected so it is locally connected using (\ref{4.6}).  Then 
$\lambda(\partial \Cal{O}^\#)=1$ and $\lambda(\partial \Cal{O}_1)=\lambda(\partial \Cal{O}_2)=0$ by Proposition \ref{p3.2}.

Because of the parabolic point at $\infty$, $R$ is neither hyperbolic nor subhyperbolic, since 
$\lim_{n \rightarrow \infty} R^n(c_1) \in \Cal{J}$ causing (\ref{4.1}) and (\ref{4.4}) to fail.

$\text{}$

\noindent{{\bf Example 2.}  Here is another example of a geometrically finite
but non hyperbolic rational map.
\beq
R(z) = \frac{z^2-z}{2z+1}.
\eeq
We list a few easily verifiable properties of $R$:
\begin{itemize}
\item  There are $3$ fixed points at $-2,0,$ and $\infty$, with respective
multipliers $1/3, -1$, and $2$.

\item  Since $-2$ is attracting and the parabolic fixed point at
$0 \in \Cal{J}$  has two attracting petals,  $\Cal{F}$ has infinitely many
components by Lemma A.1.
\item Since $R$ is quadratic, and $\Cal{J}$ is not a Cantor set,
$\Cal{J}$ is connected using A.4;  therefore each component
$\Cal{O} \subset \Cal{F}$ is simply connected.

\item Let $\Cal{O}^{\#}$ denote the component of $\Cal{F}$
containing the attracting point at $-2$.
Then since $\Cal{O}^{\#}$ contains a critical point we see that
$R:\Cal{O}^{\#} \rightarrow \Cal{O}^{\#}$  is a 2-fold branched covering of itself,
and since $R$ is degree $2$, there are no other components  of
$R^{-1} \Cal{O}^{\#}$ so it is completely invariant.

\item  The geometric finiteness comes from the fixed point at $-2$ forcing a critical point 
in $\Cal{F}$,  and similarly the fixed point at $0$ attracts the other critical point, which is also in
$\Cal{F}$.

\end{itemize}

This shows the hypotheses of Theorem \ref{t1.5} (B2) are satisfied; therefore $\lambda(\pa \Cal{O}^\#)=1$, and for any component $\Cal{O}$  mapping onto 
the components of $\Cal{F}$ containing the petals, we have $\lambda(\pa {O})=0$.

$\text{}$ \newline
{\sc Remark.}  An interesting contrast between Examples 1 and 2 lies in
the different natures of their respective completely invariant Fatou
components.

$\text{}$ \newline
{\bf Example 3.} 
We turn to a hyperbolic polynomial example.  The map $P(z)=z^2-1$ is a
hyperbolic polynomial with a superattracting period 2 cycle:
$$
0 \mapsto -1 \mapsto 0.
$$
Since $\Cal{J} = \Cal{O}^{\infty}$, $\lambda(\Cal{O}^\infty)=1$,
and there at two bounded components in the immediate basin of attraction for
the attracting 2-cycle $\{0,-1\}$.
Lemma A.1 implies that $\Cal{F}$ has infinitely many components.
By Proposition \ref{p2.3}, $\lambda(\pa \Cal{O})=0$ for all components
except $\Cal{O}^\infty$ since $\Cal{J}$ is connected and locally connected
by hyperbolicity.  There are many similar examples.

$\text{}$

\noindent {\bf Example 4.} 
There are many cubic polynomials to which Theorem \ref{1.5} (A) applies.
In particular, there are cubic polynomials with the property that one
critical point is attracted to $\infty$ and the other stays bounded.
For these we frequently see that the Julia set is disconnected but not Cantor.
Then Proposition \ref{p3.3} applies in this case.

\end{document}